\documentclass{amsart}

\usepackage{amssymb}

\usepackage[all]{xy}

\newtheorem{thm}{Theorem}

\newtheorem{lem}[thm]{Lemma}
\newtheorem{cor}[thm]{Corollary}

\newtheorem{prop}[thm]{Proposition}

\theoremstyle{definition}
\newtheorem{defn}[thm]{Definition}

\newtheorem{say}[thm]{}
\newtheorem{exmp}[thm]{Example}

\newtheorem{prob}[thm]{Problem}

\newtheorem{rem}[thm]{Remark}          

\newtheorem*{ack}{Acknowledgments}      

\newtheorem{defn-thm}[thm]{Definition--Theorem}  
\newtheorem{defn-lem}[thm]{Definition--Lemma}  

\theoremstyle{remark}

\renewcommand{\c}[0]{{\mathbb C}}  

\renewcommand{\o}[0]{{\mathcal O}} 
\newcommand{\z}[0]{{\mathbb Z}}

\renewcommand{\a}[0]{{\mathbb A}} 

\newcommand{\dd}[0]{{\mathbb D}}

\newcommand{\p}[0]{{\mathbb P}}
\newcommand{\f}[0]{{\mathbb F}}
\newcommand{\q}[0]{{\mathbb Q}}

\newcommand{\qtq}[1]{\quad\mbox{#1}\quad}
\newcommand{\spec}[0]{\operatorname{Spec}}
\newcommand{\pic}[0]{\operatorname{Pic}}

\newcommand{\supp}[0]{\operatorname{Supp}}

\newcommand{\ex}[0]{\operatorname{Ex}}

\newcommand{\depth}[0]{\operatorname{depth}} 
\newcommand{\tsum}[0]{\textstyle{\sum}}

\def\into{\DOTSB\lhook\joinrel\to}

\def\loccoh#1.#2.#3.#4.{H^{#1}_{#2}(#3,#4)}

\DeclareMathAlphabet{\mathchanc}{OT1}{pzc}%
                                {m}{it}

\newcommand{\sym}[0]{\operatorname{Sym}}

\newcommand{\defo}[0]{\operatorname{Def}}

\usepackage[all]{xy}\xyoption{dvips}

\begin{document}
\bibliographystyle{amsalpha}

\title[Grothendieck--Lefschetz type theorems]
{Grothendieck--Lefschetz type theorems \\ for the 
local Picard group}
\author{J\'anos Koll\'ar}

\maketitle

A special case of the Lefschetz hyperplane theorem asserts that
if $X$ is a smooth projective variety and $H\subset X$ an ample
divisor then the restriction map
$\pic(X)\to \pic(H)$ is an isomorphism for $\dim X\geq 4$
and an injection for $\dim X\geq 3$. 

If $X$ is normal, then the isomorphism part usually
fails. Injectivity is proved in \cite[p.305]{Kleiman66b}
and an optimal variant for the  class group  is established in
\cite{rav-sri}. 

For the local versions of these theorems,  studied in \cite{sga2},
 the projective variety is replaced by
the germ of a  singularity $(p\in X)$
and the ample
divisor by a Cartier divisor $p\in X_0\subset X$. 
The usual (global) Picard group is replaced by 
 the {\it local Picard group}
 $\pic^{\rm loc}(p\in X) $; see Definition \ref{loc.pic.defn}. 

 Grothendieck proves in \cite[XI.3.16]{sga2} that if $\depth_p \o_X\geq 4$ then
the map between the local Picard groups 
$\pic^{\rm loc}(p\in X)\to \pic^{\rm loc}(p\in X_0)$
is an injection. Note that this does not imply the Lefschetz version
since a cone over a smooth projective variety  usually has only
depth 2 at the vertex.

The aim of this note is to
 propose  a strengthening of Grothendieck's theorem that
generalizes  Kleiman's variant of the global Lefschetz theorem.
Then we prove some special cases that have interesting
applications to moduli problems.

\begin{prob} \label{main.g-l.conj}
Let $X$ be a normal (or $S_2$ and pure dimensional)
 scheme, $X_0\subset X$ a Cartier
divisor and $x\in X_0$ a closed  point. Assume that $\dim_x X\geq 4$.
What can one say about the kernel of the 
restriction  map between the local Picard groups 
$$
\operatorname{rest}^X_{X_0}:\pic^{\rm loc}(x\in X)\to \pic^{\rm loc}(x\in X_0)?
\eqno{(\ref{main.g-l.conj}.1)}
$$
We consider three  conjectural answers to
this question. 
\begin{enumerate}\setcounter{enumi}{1}
\item The map (\ref{main.g-l.conj}.1) is
 an injection if $X_0$ is $S_2$.
\item The kernel is $p$-torsion
if $X$ is an excellent, local $\f_p$-algebra.
\item The kernel  is contained in the 
connected subgroup of $\pic^{\rm loc}(x\in X) $.
\end{enumerate}
\end{prob}

The main result of this paper  gives a positive answer 
to the topological variant  (\ref{main.g-l.conj}.4) in some
cases. The precise conditions in Theorem \ref{G-L.setup.say}
are technical and they might even seem 
unrealistically special. 
Instead of stating them, I  focus on three applications first.

My interest in this subject started with 
trying to understand  higher dimensional analogs
of three theorems and examples
concerning surface singularities and their deformations;
see \cite{MR0396565, art-def, har-def} for introductions.
The main results of this note imply that none of them occurs for isolated
singularities in dimensions $\geq 3$.

\begin{say}[Three phenomena in the deformation theory  of surfaces]
\label{3.phenomena.surf}{\ }

\begin{enumerate}
\item There is a projective surface $S_0$ with quotient singularities
and ample canonical class such that  $S_0$ has
 a smoothing  $\{S_t:t\in \dd\}$ where
 $S_t$ is  a rational surface for $t\neq 0$. 

Explicit examples were written down
only recently (see \cite{lee-park, pps1, pps2} or the  simpler
\cite[3.76]{kk-singbook}), but it has been known for a long time
that $K^2$ can jump in flat families of surfaces with quotient singularities.

The simplest such example is classical and was known to Bertini
(though he probably did not consider $K^2$ for a singular surface).
Let $C_4\subset \p^5$ be the cone
over the degree 4 rational normal curve in $\p^4$. It has two
different   smoothings. In one family  $\{S_t:t\in \dd\}$
(where $S_0:=C_4$)
the general fiber is $\p^1\times \p^1\subset \p^5$ embedded by
$\o_{\p^1\times \p^1}(1,2)$. In the other family  $\{R_t:t\in \dd\}$
(where $R_0:=C_4$)
the general fiber is $\p^2\subset \p^5$ embedded by
$\o_{\p^2}(2)$.

Note that $K_{R_t}^2=9$ and $K_{S_t}^2=8$ for $t\neq 0$,
thus $K^2$ jumps in one of the families.
In fact, it is easy to compute that  $K_{C_4}^2=9$, 
so the jump occurs in the family $\{S_t:t\in \dd\}$.

\item There are non-normal, isolated, smoothable  surface singularities
$(0\in S_0)$ whose normalization is simple elliptic \cite{MR541025}.

\item  Every rational surface singularity
$(0\in S_0)$ has a smoothing that admits a simultaneous resolution.

It is known that such smoothings form a whole component of the deformation space
$\defo(0\in S_0)$, called the {\it Artin component}
\cite{art-bri}.  Generalizations of this can be used to describe
all components of the deformation space of quotient singularities
\cite{ksb}, and, conjecturally, of all rational surface singularities
and many other non-rational singularities \cite{Kollar91, MR1163728}.
\end{enumerate}
\end{say}

The higher dimensional versions of these were studied with
the ultimate aim of compactifying the moduli space of
varieties of general type; see \cite{k-modsurv} for an introduction.
The general theory of \cite{ksb, k-modsurv}
suggests that one should work with
 {\it log canonical} singularities; see \cite{km-book} or
(\ref{lc.defn}) for their definition and  basic properties.
This guides our generalizations of (\ref{3.phenomena.surf}.1--2).

In order to develop  (\ref{3.phenomena.surf}.3) further, recall  that a 
surface singularity is rational iff its divisor class group is finite; see 
\cite{mumf-top}.

In this paper  we state the results for normal varieties.
In the  applications to moduli questions
one needs these results for semi-log canonical pairs
$(X,\Delta)$. All  the theorems extend to this general setting
using the methods of \cite[Chap.5]{kk-singbook};
see the forthcoming \cite[Chap.3]{k-modbook} for details.

We say that a variety $W$ is {\it smooth} (resp.\ {\it normal})
{\it in codimension $r$} if there is a closed subscheme $Z\subset W$
of codimension $\geq r+1$ such that  $W\setminus Z$ is  smooth
 (resp.\  normal).

\begin{thm} \label{3.phenomena.hd}
None of the above examples (\ref{3.phenomena.surf}.1--3) exists for
varieties with isolated singularities in dimension $\geq 3$. 
More generally  the following hold.
\begin{enumerate}
\item  Let $X_0$ be a projective variety  
with  log canonical singularities
and ample canonical class. 
If $X_0$ is smooth in codimension 2 then
every deformation of $X_0$ also has ample canonical class.
\item Let $X_0$ be a  non-normal variety
 whose normalization is log canonical. 
If $X_0$ is normal in codimension 2
then $X_0$ is not smoothable,
it does not even have normal deformations.

\item  Let $ X_0$ be a normal variety
whose local class groups are torsion and
 $\{X_t:t\in \dd\}$ a smoothing. 
If $X_0$ is smooth in codimension 2 then
 $\{X_t:t\in \dd\}$ does not admit a simultaneous resolution.
\end{enumerate}
\end{thm}

Our results in Sections \ref{appl.sec}--\ref{exdivdef.sec}
 are even stronger; we need only some
control over the singularities in codimension 2.

Such results have been known if $X_0$ (or its normalization) has
rational singularities. These essentially follow from \cite[XI.3.16]{sga2};
see \cite{k-flat} for details. Thus the new part of Theorem 
\ref{3.phenomena.hd}.1--2
is that the claims also hold for log canonical singularities
that are not rational.

\subsection*{Further comments and problems}\footnote{Recently Bhatt and de~Jong proved Conjecture \ref{main.g-l.conj}.3 in general
and Conjecture \ref{main.g-l.conj}.2 for schemes essentially of finite type in
characteristic 0.}{\ }

The theorems of SGA  rarely have unnecessary assumptions,
so an explanation is needed why  Problem \ref{main.g-l.conj}
could be an exception. One reason is that while our  assumptions
are weaker, the conclusions in  \cite{sga2} 
are stronger.

\begin{thm} \cite[XI.2.2]{sga2}
\label{fibw.cartier.codim3.G.thm} 
Let $X$ be a scheme of pure dimension $n+1$,
$X_0\subset X$ a Cartier divisor and
$Z\subset X$ a closed subscheme such that
$Z_0:=X_0\cap Z$ has  dimension $\leq n-3$.
Let $D^*$ be  a Cartier divisor on $X\setminus Z $
such that $D^*|_{X_0\setminus Z_0}$ extends to a Cartier divisor on $X_0$.
Assume furthermore that $X_0$ is $S_2$ and $\depth_{Z_0}\o_{X_0}\geq 3$.
Then $D^*$ extends to a Cartier divisor on $X$, in some neighborhood of $X_0$.
\qed
\end{thm}

 In Problem \ref{main.g-l.conj} we assume that 
$Z$ is contained in $X_0$, thus it is not entirely
surprising that the depth condition $\depth_{Z_0} X_0\geq 3$,
could be relaxed.

  Example \ref{no.str.GL.forcones} shows that in
Theorem \ref{fibw.cartier.codim3.G.thm}  the condition $\depth_{Z_0} X_0\geq 3$
is necessary.

Another reason why Problem \ref{main.g-l.conj} may have escaped
attention is that the topological version of it fails.
 In (\ref{Lef.fails.for.H2.lem}--\ref{Lef.fails.for.H2.exmp}) we
 construct normal, projective varieties
$Y$ (of arbitrary large dimension) with a single singular point $y\in Y$ and
a smooth hyperplane section $H$ (not passing through $y$)
such that the restriction map
$H^2(Y, \q)\to H^2(H, \q)$ is not injective.
However, the kernel does not contain $(1,1)$-classes. 
In the example $Y$ even has a log canonical singularity at $y$.

The arguments in Section \ref{top-approach.sec}
 show that, at least over $\c$,  
a solution to Problem \ref{main.g-l.conj}
would be implied by the following.

\begin{prob} \label{main.g-l.stein.conj}
Let $W$ be a normal Stein space of dimension $\geq 3$
and $L$ a holomorphic line bundle on $W$. Assume that 
there is a compact set $K\subset W$ such that
the restriction of $c_1(L)$ is zero in $H^2(W\setminus K, \z)$.

Does this imply that $c_1(L)$ is zero in $H^2(W, \z)$?
\end{prob}

Another approach would be to use intersection cohomology
to restore Poincar\'e duality in (\ref{leff.in.intersection.hom}.2). 
For this to work, the solution of the following is needed.

\begin{prob} \label{main.IH.prob}
Let $W$ be a normal analytic space. Is there an exact sequence
$$
H^1(W, \o_W)\to \pic^{\rm an}(W)\otimes\q\to IH^2(W, \q)?
\eqno{(\ref{main.IH.prob}.1)}
$$
\end{prob}

Note first that the sequence exists.
 Indeed, let $g:W'\to W$
be a resolution of singularities. If $L$ is a line bundle on $W$,
then $g^*L$ is a line bundle on $W'$ hence it has a Chern class
$c_1\bigl(g^*L)\in H^2(W', \z)$. By the decomposition theorem
\cite{MR751966}
$IH^2(W, \q)$ is a direct summand of $ H^2(W', \q)$
(at least for algebraic varieties).

Arapura explained to me that the sequence (\ref{main.IH.prob}.1)
should be exact for projective varieties by weight considerations
but the general complex case is not clear.
For our applications we need the case when $W$ is Stein.

 \begin{ack}
I thank D.~Arapura,  C.~Xu and R.~Zong
 for discussions and remarks.
 The formulation of Conjectures \ref{main.g-l.conj}.2--4
owes a lot to the comments and results of 
B.~Bhatt and  A.~J.~de~Jong.
The referees comments helped to eliminate several inaccuracies.
Partial financial support   was provided  by  the NSF under grant number 
DMS-07-58275 and by the Simons Foundation. Part of this paper was written 
while the author visited the University of Utah.
\end{ack}

\section{Definitions and examples}

\begin{defn}[Local Picard groups]\label{loc.pic.defn}
Let $X$ be a scheme and $p\in X$ a point. 
The  {\it local Picard group}
 $\pic^{\rm loc}(p\in X) $ is a group whose 
 elements are $S_2$ sheaves $F$ on some neighborhood
$p\in U\subset X$ such that $F$ is locally free on $U\setminus \{p\}$.
Two such sheaves give the same element if they are isomorphic over some
 neighborhood of $p$. The product is given by the $S_2$-hull of the
tensor product.

One can also realize the local Picard group as the direct limit of
 $\pic(U\setminus\{p\})$ as $U$ runs through all open 
Zariski neighborhoods of $p$ or as
$\pic(\spec \o_{x,X}\setminus\{p\})$. 

If $X$ is normal and $X\setminus \{x\}$ is smooth then  $\pic^{\rm loc}(p\in X) $
is isomorphic to  the divisor class group of $\o_{x,X}$.

In many contexts it is more natural to work with the 
{\it \'etale-local Picard group} $\pic^{\rm et-loc}(p\in X) :=
\pic(\spec \o_{x,X}^h\setminus\{p\})$ where $\o_{x,X}^h $
is the Henselization of the local ring $\o_{x,X} $.

 If $X$ is defined over $\c$, let $W\subset X$ be the intersection 
of $X$ with a small (open) ball around $p$. The 
{\it analytic local Picard group} $\pic^{\rm an-loc}(p\in X)$
can be defined as above using (analytic) $S_2$ sheaves on $W$. 
By \cite{Artin69d}, there is a natural isomorphism
$$
\pic^{\rm et-loc}(p\in X)\cong \pic^{\rm an-loc}(p\in X).
$$
Note that $\pic^{\rm an}(W\setminus\{p\})$ is usually
much bigger than $\pic^{\rm an-loc}(p\in W) $. (This happens already for
$X=\c^2$.) However, 
$\pic^{\rm an}(W\setminus\{p\})=\pic^{\rm an-loc}(p\in X) $ if $\depth_p\o_X\geq 3$
 \cite{MR0245837}.

(The literature does not seem to be consistent; any of the above
four variants is called the  local Picard group by some authors.)
\end{defn}

Let $X$ be  a complex space and $p\in X$ a closed point.
Set $U:=X\setminus \{p\}$. As usual, 
$\pic(U)\cong H^1(U, \o_U^*)$ and the exponential sequence
$$
0\to \z_U\stackrel{2\pi i}{\longrightarrow} \o_U
\stackrel{\rm exp}{\longrightarrow} \o_U^*\to 1
$$
gives  an exact sequence
$$
H^1\bigl(U, \o_U\bigr)\to \pic(U) \stackrel{c_1}{\longrightarrow} 
H^2\bigl(U, \z\bigr).
$$
A piece of  the  local cohomology  exact sequence is
$$
 H^1\bigl(X, \o_X\bigr)\to 
H^1\bigl(U, \o_U\bigr)\to H^2_p\bigl(X, \o_X\bigr)\to H^2\bigl(X, \o_X\bigr).
$$
Thus if  $X$ is Stein  then we have an isomorphism
$$
H^1\bigl(U, \o_U\bigr)\cong H^2_p\bigl(X, \o_X\bigr)
$$
and the latter vanishes iff $\depth_p\o_X\geq 3$; see
 \cite[Sec.3]{gro-loc-coh-MR0224620}. 
Combining with  \cite{MR0245837} we obtain the following well known
result.

\begin{lem} \label{pic.H2.lem}
Let $X$ be a  Stein space and $p\in X$
 a closed point. If $\depth_p\o_X\geq 3$
then taking the first Chern class gives an injection
$$
c_1: \pic^{\rm an-loc}(p\in X)\into 
H^2\bigl(X\setminus \{p\}, \z\bigr).
\qed
$$
\end{lem}

\begin{defn}[Log canonical singularities]\label{lc.defn}
(See \cite{km-book} for an introduction and \cite{kk-singbook}
for a comprehensive treatment of these singularities.)

Let $X$ be a normal variety such that $mK_X$ is Cartier for some $m>0$.
Let $f:Y\to X$ be a resolution of singularities
with exceptional divisors $\{E_i:i\in I\}$. One can then write
$$
mK_Y\sim f^*(mK_X)+m\cdot \tsum_{i\in I} a(E_i, X)E_i.
$$
The number $a(E_i, X)\in \q$ is called the
{\it discrepancy} of $E_i$; it is independent of the choice of $m$.

If $\min\{a(E_i, X):i\in I\}\geq -1$ then the minimum is
independent of the resolution $f:Y\to X$ and its value is called
the {\it discrepancy} of $X$.

$X$ is called {\it log canonical}  if $\min\{a(E_i, X):i\in I\}\geq -1$
and {\it log terminal}  if $\min\{a(E_i, X):i\in I\}> -1$.
A cone over a smooth variety with trivial canonical class is
log canonical but not log terminal.

Let $X$ be log canonical, $g:X'\to X$ any resolution
and $E\subset X'$ an exceptional divisor
such that $a(E, X)=-1$. The 
 subvariety $g(E)\subset X$ is called a {\it log canonical center} of $X$. 
Log canonical centers hold the key to understanding
log canonical varieties, see \cite[Chaps.4--5]{kk-singbook}.

Log terminal singularities are rational \cite[5.22]{km-book}.
Log canonical  singularities are usually not rational but they
 are Du~Bois  \cite{k-db}.

Log canonical singularities
(and their non-normal versions, called semi-log canonical singularities)
are precisely those that are needed to compactify the
moduli of varieties of general type.
\end{defn}

We use only the following two theorems about log canonical singularities.

\begin{thm}\label{inv.adj.thms}\cite{MR2264806, ale-lim}
Let $X$ be a normal variety over a field of characteristic 0,
$Z\subset X$ a closed subscheme of codimension $\geq 3$ and 
 $Z\subset X_0\subset X$ a Cartier divisor 
such that $X_0\setminus Z$ is normal and the normalization of $X_0$ is
  log canonical.
Assume also that $K_X$ is $\q$-Cartier. Then $X_0$ is normal
and it does not contain any log canonical center of $X$. \qed
\end{thm}

\begin{thm}\label{small.res.thms}
Let $X$ be a normal variety over a field of characteristic 0,
$Z\subset X$ a closed subscheme of codimension $\geq 3$
and $D$ a Weil divisor on $X$ that is Cartier on $X\setminus Z$.

Then there is a proper, birational morphism $f:Y\to X$
such that $f$ is small (that is, its exceptional set has codimension $\geq 2$)
and the birational transform $f^{-1}_*D$ is $\q$-Cartier and $f$-ample
if one of the following assumptions is satisfied.
\begin{enumerate}
\item \cite{ksb} There is a Cartier divisor  $Z\subset X_0\subset X$
such that $X_0\setminus Z$  is normal and the normalization of $X_0$ is
   log canonical.
\item \cite{birkar11, hacon-xu, oda-xu} $X$ is log canonical
and $Z$ does not contain any log canonical center of $X$. \qed
\end{enumerate}
\end{thm}

Observe that $D$ is  $\q$-Cartier iff $f:Y\to X$ is an isomorphism.
In our applications we show that $Y\neq X$ leads to a contradiction.
\medskip

{\it Comments on the references.} \cite{MR2264806} proves that
the pair  $(X, X_0)$ is log canonical. This implies that 
 $X_0$  does not contain any log canonical center of $X$
by an easy monotonicity argument  \cite[2.27]{km-book}.
Then \cite{ale-lim} shows that $X_0$ is $S_2$, hence normal since
we assumed normality in codimension 1.

 \cite{ksb} claimed (\ref{small.res.thms}.1) only for $\dim X=3$ since
the necessary results of Mori's program were known only 
for $\dim X\leq 3$ at that time. The proof of the general case is the same.

The second case (\ref{small.res.thms}.2)
 is not explicitly claimed in the references
but it easily follows from them. For details on both cases
see \cite[Sec.1.4]{kk-singbook}.
\medskip

  The next example  shows that 
Theorem \ref{fibw.cartier.codim3.G.thm}  fails if  $\depth_{Z_0} X_0<3$,
even if the dimension is large.

\begin{exmp}  \label{no.str.GL.forcones} 
Let $(A,\Theta)$ be a principally 
polarized Abelian variety over a field $k$.
The  affine cone over $A$
with vertex $v$ is
$$
C_a(A,\Theta):=\spec_k \tsum_{m\geq 0} H^0\bigl(A, \o_A(m\Theta)\bigr).
$$ 
Note that $\depth_{v} C_a(A,\Theta)=2$ since
$H^1(A, \o_A)\neq 0$.

Set $X:=C_a(A,\Theta)\times \pic^0(A)$ with
$f:X\to \pic^0(A)$ the second projection.
Since $L(\Theta)$ has a unique section for every $L\in  \pic^0(A)$,
there is a unique divisor $D_A$ on $A\times  \pic^0(A)$
whose restriction to $A\times \{[L]\}$ is the above divisor.
By taking the cone we get a divisor $D_X$ on $X$.

For $L\in  \pic^0(A)$, let $D_{[L]}$ denote the restriction of $D_X$
to the fiber $C_a(A,\Theta)\times \{[L]\}$ of $f$.
We see that
\begin{enumerate}
\item $D_{[L]}$  is Cartier iff $L\cong \o_A$.
\item $mD_{[L]}$  is Cartier iff $L^m\cong \o_A$.
\item $D_{[L]}$  is not $\q$-Cartier for very general $L\in  \pic^0(A)$.
\end{enumerate}
\end{exmp}


\section{The main technical theorem}

The following is our main result concerning Problem \ref{main.g-l.conj}.
In the applications the key question will be the existence of the 
bimeromorphic morphism $f:Y\to X$.
This is a very hard question in general, but in our cases 
existence is guaranteed by  Theorem \ref{small.res.thms}.

\begin{thm} \label{G-L.setup.say}
Let $f:Y\to X$ be a proper, bimeromorphic morphism of 
normal analytic spaces of dimension $\geq 4$ and $L$ 
 a line bundle on $Y$ whose restriction to every fiber  is ample.
Assume that there is  a closed subvariety $Z_Y\subset Y$
 of codimension $\geq 2$ such that 
$Z:=f\bigl(Z_Y\bigr)$ has dimension $\leq 1$ and
$f$ induces an isomorphism $Y\setminus Z_Y\cong X\setminus Z$.

Let $X_0\subset X$
be a Cartier divisor such that $Z\cap X_0$ is a single point $p$.
Let  $p\in U\subset X$ be a 
contractible open neighborhood of $p$. 
Note that $U\setminus Z\cong f^{-1}(U)\setminus Z_Y$, hence
the restriction $L|_{U\setminus Z}$ makes sense.
Set $U_0:=X_0\cap U$. 
The following are equivalent.
\begin{enumerate}
\item The map $f$ is an isomorphism over $U$.
\item  The Chern class of $L|_{U\setminus Z}$ vanishes in
$H^2\bigl( U\setminus Z, \q\bigr)$.
\item  The Chern class of $L|_{U_0\setminus \{p\}}$ vanishes in
$H^2\bigl( U_0\setminus \{p\}, \q\bigr)$.
\end{enumerate}
\end{thm}

Proof. If $f$ is an isomorphism then $L$ is a line bundle
on the contractible space $U$ hence $c_1(L)=0$ in
$H^2\bigl( U, \q\bigr)$. Thus (2) holds and clearly (2) implies (3). 
The key part is to prove that (3) $\Rightarrow$ (1).

The assumption and the conclusion are both local near $p$
in the Euclidean topology. By shrinking $X$ we may assume that
the Cartier divisor $X_0$ gives a morphism $g:X\to \dd$ 
to the unit disc $\dd$ whose central
fiber is $X_0$. Note that $Y_0:=f^{-1}(X_0)$
is a Cartier divisor in $Y$.

Let $W\subset X$ be the intersection of
$X$ with a  closed ball of radius $0<\epsilon \ll 1$ around $p$.
Set $W_0:=X_0\cap W$,   $V:=f^{-1}(W)$ and $V_0:=f^{-1}(W_0)$.
We may assume that $W, W_0$ are contractible and
$f^{-1}(p)$ is a strong deformation retract of both $V$ and of $V_0$.

Let $\bar \dd_{\delta}\subset \dd $ denote the 
closed disc of radius $\delta$.
If $0< \delta\ll \epsilon$ then  the pair
$\bigl( W_0, \partial  W_0\bigr)$ is a  strong deformation retract of
$\bigl( W\cap g^{-1}\bar \dd_{\delta},   
\partial W\cap g^{-1}\bar \dd_{\delta}\bigr)$
and
$\bigl( V_0, \partial  V_0\bigr)$
is a  strong deformation retract of
$\bigl( V\cap (gf)^{-1}\bar \dd_{\delta},   
\partial V\cap (gf)^{-1}\bar \dd_{\delta}\bigr)$.
These retractions induce continuous maps
(unique up-to homotopy)
$$
r_c: \bigl( V_c, \partial  V_c\bigr) \to \bigl( V_0, \partial  V_0\bigr),
\eqno{(\ref{G-L.setup.say}.4)}
$$
where $V_c$ is the fiber of $g\circ f:V\to \dd$ over $c\in \dd_{\delta}$.
The induced maps
$$
r_c^*: \q\cong H^{2n}\bigl( V_0, \partial  V_0, \q\bigr)
\stackrel{\cong}{\to}
H^{2n}\bigl( V_c, \partial  V_c,\q, \bigr)\cong \q
\eqno{(\ref{G-L.setup.say}.5)}
$$
are isomorphisms where $n=\dim_{\c} V_0=\dim_{\c} V_c$.
Our aim is to study the 
cup product pairing
$$
H^{2}\bigl( V_0, \partial  V_0, \q\bigr)\times
H^{2n-2}\bigl( V_0,  \q\bigr)\to 
H^{2n}\bigl( V_0, \partial  V_0, \q\bigr)\cong \q.
\eqno{(\ref{G-L.setup.say}.6)}
$$
(See \cite{Hatcher}, especially pages 209 and 240 for the
relevant facts on cup and cap products.)
We prove in Lemmas \ref{G-L.setup.vanish.lem}
and \ref{G-L.setup.vanish.lem2}, by arguing on $V_c$, that it is zero 
and in Lemma \ref{G-L.setup.nonvanish.lem}.2, by arguing on $V_0$, that it is
 nonzero  if $V_0\to W_0$ is not finite.
Thus $f^{-1}(p)$ is 0-dimensional, hence $f$ is a biholomorphism. 
\qed
\medskip

For  later applications,
in the next lemmas we consider the more general case when
 $f:Y\to X$ is a proper, bimeromorphic morphism of normal analytic spaces
and $g:X\to \dd$  a flat morphism 
of relative dimension $n$.

\begin{lem} \label{G-L.setup.vanish.lem}
Notation and assumptions as in 
{\rm (\ref{G-L.setup.say})}.  If $H_{2n-2}\bigl( V_c,\q \bigr)=0$
then the cup product pairing
$$
H^{2}\bigl( V_0, \partial  V_0, \q\bigr)\times
H^{2n-2}\bigl( V_0,  \q\bigr)\to 
H^{2n}\bigl( V_0, \partial  V_0, \q\bigr)\cong \q
\qtq{is zero.}
$$
\end{lem}

Proof. Using $r_c^*$ and the Poincar\'e duality map,
the cup product pairing factors through the following
cup and cap product pairings, where the 
right hand sides are isomorphic by 
(\ref{G-L.setup.say}.5),
$$
\begin{array}{ccccccc}
H^{2}\bigl( V_0, \partial  V_0, \q\bigr)&\times&
H^{2n-2}\bigl( V_0,  \q\bigr)&\to &
H^{2n}\bigl( V_0, \partial  V_0, \q\bigr)& \cong & \q\\
\downarrow && \downarrow &&\hphantom{\cong}\downarrow & &||  \\
H^{2}\bigl( V_c, \partial  V_c, \q\bigr)&\times &
H^{2n-2}\bigl( V_c,  \q\bigr)&\to & 
H^{2n}\bigl( V_c, \partial  V_c, \q\bigr)& \cong & \q\\
\downarrow && \downarrow &&\hphantom{\cong}\downarrow & &||  \\
H_{2n-2}\bigl( V_c, \q \bigr)&\times &
H^{2n-2}\bigl( V_c,  \q\bigr)&\to & H_{0}\bigl( V_c,  \q\bigr)& \cong & \q
\end{array}
$$
The first factor in the bottom row is zero, hence the
pairing is zero. \qed
\medskip

We apply the next result to  $V_c\to W_c$ to check the
homology vanishing assumption in Lemma \ref{G-L.setup.vanish.lem}.

\begin{lem} \label{G-L.setup.vanish.lem2}
Let $V'\to W'$ be a proper bimeromorphic map of normal complex spaces
of dimension $n\geq 3$.
Assume that  every fiber  has complex dimension
$\leq n-2$ and $W'$ is Stein. 
Then $H_{2n-2}\bigl( V', \q \bigr)=0$.
\end{lem}

Proof.
Let $E'\subset V'$ denote the exceptional set
and $F'\subset W'$ its image. Then $\dim F'\leq n-2$, hence
the exact sequence
$$
H_{2n-2}\bigl( F', \q \bigr)\to
H_{2n-2}\bigl( W', \q \bigr)\to H_{2n-2}\bigl( W', F', \q \bigr)
\to H_{2n-3}\bigl( F', \q \bigr)
$$
shows that 
$ H_{2n-2}\bigl( W', \q \bigr)\cong H_{2n-2}\bigl( W', F', \q \bigr)$.
The latter group is in turn isomorphic to
$H_{2n-2}\bigl( V', E', \q \bigr) $
which sits in an exact sequence
$$
H_{2n-2}\bigl( E', \q \bigr)\to
H_{2n-2}\bigl( V', \q \bigr)\to H_{2n-2}\bigl( V', E', \q \bigr).
$$
Here $H_{2n-2}\bigl( E', \q \bigr) $
is generated by the fundamental classes of the compact irreducible
components of $E'$, but we assumed that there are no such.
Thus we have an injection
$$
H_{2n-2}\bigl( V', \q \bigr)\into H_{2n-2}\bigl( W', \q \bigr).
$$
 Since $W'$ is Stein and  $2n-2>n$,  Theorem \ref{stein.top.main.thm}
 implies that 
$H_{2n-2}\bigl( W', \q \bigr) =0$.
Thus we conclude that $H_{2n-2}\bigl( V', \q \bigr) =0$. \qed
\medskip

During the proof we have used the following.

\begin{thm} \label{stein.top.main.thm} \cite{MR684017, MR817633} or
\cite[p.152]{gm-book}.
Let $W$ be a Stein space of dimension $n$. Then $H_i(W,\z)$ and $H^i(W,\z)$
both vanish for $i>n$. 
More generally, $W$ is homotopic to a CW complex of dimension $\leq n$.\qed
\end{thm}

Next we describe two cases when the cup product pairing
(\ref{G-L.setup.say}.6) is nonzero.
The first of these is used in  Proposition  \ref{stab.exc.divs.cor}
and the second in  Theorem  \ref{G-L.setup.say}. 

\begin{lem} \label{G-L.setup.nonvanish.lem}
 Let $f_0:V_0\to W_0$ be a projective,
bimeromorphic morphism
between irreducible complex spaces. Let $p\in W_0$ be a point. Assume that
$f_0$ is an isomorphism over $W_0\setminus \{p\}$
and $\dim f_0^{-1}(p)>0$.
Assume furthermore that 
one of the following holds.
\begin{enumerate}
\item There is a nonzero
$\q$-Cartier divisor $E_0\subset V_0$ supported on
$f_0^{-1}(p)$.  
\item  There is an $f_0$-ample
line bundle $L$ such that
$c_1(L)|_{\partial V_0}=0$.
\end{enumerate}
Then  $f_0^{-1}(p)$ has codimension 1 in $ V_0$ and
 the cup product pairing
$$
H^{2}\bigl( V_0, \partial  V_0, \q\bigr)\times
H^{2n-2}\bigl( V_0,  \q\bigr)\to 
H^{2n}\bigl( V_0, \partial  V_0, \q\bigr)\cong \q
\qtq{is nonzero.}
$$
\end{lem}

Proof. Consider first Case (1). Then $E_0\neq 0$ shows that
$f_0^{-1}(p)$ has codimension 1.

Let $H$ be a relatively very ample  line bundle.
We have  $c_1(E_0)\in H^{2}\bigl( V_0, \partial  V_0, \q\bigr)$
and $ c_1(H)\in H^{2}\bigl( V_0,  \q\bigr)$.
If $E_0$ is effective then
$$
c_1(E_0)\cup c_1(H)^{n-1}= c_1\bigl(H|_{E_0}\bigr)^{n-1}
\in H_0\bigl( E_0, \q\bigr)\to H_0\bigl( V_0, \q\bigr)
\eqno{(\ref{G-L.setup.nonvanish.lem}.3)}
$$
is positive. If $E_0$ is not assumed effective then
we claim that
$$
c_1(E_0)\cup c_1(E_0)\cup 
c_1(H)^{n-2}\in H^{2n}\bigl( V_0, \partial  V_0, \q\bigr)
\qtq{is nonzero.}
$$
The complete intersection of $(n-2)$ general members of $H$
gives an algebraic surface $S$, proper over $W_0$
such that $E_0\cap S$ is a nonzero linear combination of 
exceptional curves. Thus, by the Hodge index theorem,
$$
c_1(E_0)^2\cup c_1(H)^{n-2}=c_1\bigl(E_0|_S\bigr)^2 <0,
$$
completing the proof of (1). 

Next assume that (2) holds.
 By assumption we can lift $c_1(L)$ to
$\tilde c_1(L)\in H^{2}\bigl( V_0, \partial  V_0, \q\bigr)$. 
(The lifting is in fact unique, but this is not important for us.)
From this we obtain a class
$$
\bigl[\tilde c_1(L)\bigr]\in H_{2n-2}\bigl( V_0, \q\bigr)=
H_{2n-2}\bigl( f^{-1}(p), \q\bigr)=\tsum \q[A_i]
\eqno{(\ref{G-L.setup.nonvanish.lem}.4)}
$$
where $A_i\subset f^{-1}(p)$ are the irreducible components of
dimension $n-1$. So far we have not established that
$\dim f^{-1}(p)=n-1$, thus the sum in (\ref{G-L.setup.nonvanish.lem}.4) 
could be empty.
The key step is the following.
\medskip

{\it Claim} \ref{G-L.setup.nonvanish.lem}.5.
$\bigl[\tilde c_1(L)\bigr]=\sum a_i[A_i]$
where  $a_i<0$ for every $i$ and the sum is not empty. 
\medskip

Once this is shown
we conclude as in  (\ref{G-L.setup.nonvanish.lem}.3) using the equality
$$
\tilde c_1(L)\cup c_1(L)^{n-1}= \tsum a_i \cdot c_1\bigl(L|_{A_i}\bigr)^{n-1}<0.
\eqno{(\ref{G-L.setup.nonvanish.lem}.6)}
$$

In order to prove (\ref{G-L.setup.nonvanish.lem}.5)
we aim to use \cite[3.39]{km-book}, except there it is assumed that
every $A_i$ is $\q$-Cartier.  To overcome this, take a resolution
$\pi: V'_0\to V_0$ and write the 
homology class $\bigl[\tilde c_1(\pi^*L)\bigr]$
 is a linear combination
$\sum a'_i[A'_i]$ where 
$A'_i\subset (f\circ \pi)^{-1}(p)$ are the irreducible components of
dimension $n-1$.
Since $L$ is $f$-ample and $\dim f^{-1}(p)>0$,
we see that $\pi^*L$ is nef and not numerically trivial on 
$(f\circ \pi)^{-1}(p)$.
Apply \cite[3.39.2]{km-book} to $\pi^*L$.
We obtain that  $-\sum a'_i[A'_i]$ is effective and its
support contains  $(f\circ \pi)^{-1}(p)$. Thus 
$a'_i<0$ for every $i$
and so  $\bigl[\tilde c_1(L)\bigr]=\pi_*\sum a'_i[A'_i]$
shows (\ref{G-L.setup.nonvanish.lem}.5) unless there are
no $f$-exceptional divisors $A_i\subset f^{-1}(p)$.

If this happens, then $\sum a'_i[A'_i]$ is  $\pi$-exceptional
and, as the homology class of   $\pi^*L$, it  has zero intersection
with   every curve that is
contracted by $\pi$. Thus we can apply  \cite[3.39.1]{km-book}
to both $\pm \sum a'_i[A'_i]$ and conclude that
$\sum a'_i[A'_i]=0$. This is a contradiction since
$L$ and hence $\pi^*L$ have positive intersection with some
curve. \qed
\medskip

\section{Deformations of log canonical singularities}
\label{appl.sec}

Here we derive stronger forms of the three assertions
of Theorem \ref{3.phenomena.hd}.
We start with (\ref{3.phenomena.hd}.1--2).

\begin{thm}\label{codim.3.slc.GL.thm} Let $X$ be a normal variety
over $\c$ 
 and $g:X\to C$  a flat morphism
of pure relative dimension $n$  to a smooth curve.
 Let $0\in C$ be a point and
$Z_0\subset X_0$ a closed subscheme of dimension $\leq n-3$.
Assume that 
\begin{enumerate}
\item $K_X$ is  $\q$-Cartier on $X\setminus Z_0$,
\item the fibers $X_c$ are log canonical for $c\neq 0$,
\item  $X_0\setminus Z_0$ is log canonical and
\item the normalization of $X_0$ is log canonical.
\end{enumerate}
Then  $X_0$ is  normal and $K_X$ is  $\q$-Cartier on $X$.
\end{thm}

Proof. 
 By localization we may assume
that $Z_0=\{p\}$ is closed point.
Next we use Theorem \ref{small.res.thms} to obtain
$f:Y\to X$ such that $f$ is an isomorphism over
$X\setminus \{p\}$ and $f^{-1}_*K_X$ is
an $f$-ample $\q$-Cartier divisor.

By the Lefschetz principle, we may assume that everything is defined
over $\c$.

We apply Theorem \ref{G-L.setup.say}
 to $L:=mf^{-1}_*K_X$ for a suitable $m>0$.
Let $U$ be the intersection of $X$ with a small ball around $p$
and set $U_0:=X_0\cap U$. Note that $U_0$ is naturally a subset
of $X$, of $Y$ and also of the normalization  of $ X_0$.
The latter shows that $L_{U_0\setminus\{p\}}$ is trivial,
thus the assumption (\ref{G-L.setup.say}.3) is satisfied. 
Hence $f$ is an isomorphism
and so $K_X$ is $\q$-Cartier.

Now Theorem \ref{inv.adj.thms} implies that $X_0$ is normal. \qed
\medskip

\begin{thm}\label{GL.pic-inj.nonlc.thm}
 Let $X$ be a log canonical variety of dimension $\geq 4$ over $\c$
and $p\in X$ a closed point that is not a log canonical center (\ref{lc.defn}). 
Let  $p\in X_0\subset X$ be  a Cartier divisor.
Let $p\in U$ be a Stein neighborhood such that
$U$ and $U_0:=X_0\cap U$ are both contractible.
 Then the restriction maps
 $$
\pic^{\rm loc}(p\in X)\to \pic^{\rm loc}(p\in X_0)
\qtq{and} \pic^{\rm loc}(p\in X)\to H^2\bigl(U_0\setminus \{p\}, \z\bigr)
$$
are injective.  
\end{thm}

Proof. Let 
$D$ be a divisor on $X$ such that $D|_{X\setminus \{p\}}$ is Cartier
 and  $c_1\bigl( D|_{X_0\setminus \{p\}}\bigr)$ is zero in 
$H^2\bigl(X_0\setminus \{p\}, \z\bigr)$.

First we show that $D$ is $\q$-Cartier at $p$.
By Theorem \ref{small.res.thms}
there is a proper birational morphism
$f:Y\to X$ such that $f$ is an isomorphism over  $X\setminus \{p\}$,
$f^{-1}_*D$ is $f$-ample and $f$ has no exceptional divisors.

The Cartier divisor $X_0$ gives a morphism $X\to \dd$ whose central
fiber is $X_0$. 

As in Theorem \ref{G-L.setup.say} let  $W_0\subset X_0$ be the intersection of
$X_0$ with a  closed ball of radius $\epsilon$ around $p$ and $V_0:=f^{-1}(W_0)$.
Set $n:=\dim X_0$.
By Lemma \ref{G-L.setup.vanish.lem} we see that the cup product pairing
$$
H^{2}\bigl( V_0, \partial  V_0, \q\bigr)\times
H^{2n-2}\bigl( V_0,  \q\bigr)\to 
H^{2n}\bigl( V_0, \partial  V_0, \q\bigr)\cong \q
\qtq{is zero. }
$$
 On the other hand, by (\ref{G-L.setup.nonvanish.lem}.2)
 it is nonzero unless
 $f:Y\to X$ is finite. Thus $f$ is
 an isomorphism and $D$ is $\q$-Cartier at $p$. 

Now we can use \cite[X.3.2]{sga2}  
 to show that   $D$ is Cartier at $p$.
\qed
\medskip

\begin{cor}\label{fibw.cartier.codim3.thm} 
Let $g:X\to C$  and $Z_0\subset X_0$ be as in  Theorem \ref{codim.3.slc.GL.thm}.
Assume that the fibers $X_c$ are all log canonical and
$K_X$ is $\q$-Cartier.
Let $D^*$ be  Cartier divisor on $X\setminus Z_0 $
such that $D^*|_{X_0\setminus Z_0}$ extends to a Cartier divisor on $X_0$.

Then $D^*$ extends to a Cartier divisor on $X$.
\end{cor}

Proof. Choose $Z$ to be the smallest closed
subset such that $D^*$ is Cartier on $X\setminus Z$. We need to show that
$Z=\emptyset$. If not, let $p\in Z$ be a generic point. By localization
we are reduced to the case when $Z=\{p\}$ is a closed point of $X$. 
Note that $p$ is not a log canonical
 center of $X$  by Theorem \ref{inv.adj.thms}. 

Thus (\ref{fibw.cartier.codim3.thm})
is a special case of Theorem \ref{GL.pic-inj.nonlc.thm}. 
\qed

\begin{say}[Proof of Theorem \ref{3.phenomena.hd}.1--2]
Let
 $g:X\to C$ be a flat, proper morphism to a smooth curve whose fibers are
normal and log canonical. Let $0\in C$ be a closed  point and  $Z_0\subset X_0$ 
a subscheme of 
codimension $\geq 3$ such that $K_X$ is $\q$-Cartier on $X\setminus Z_0$.
Then $K_X$ is $\q$-Cartier by Corollary \ref{fibw.cartier.codim3.thm}
thus $mK_X$ is Cartier  for some $m>0$.
So $\o_X(mK_X)$ is a line bundle on $X$.
For a flat family of line bundles, ampleness is an open condition,
proving (\ref{3.phenomena.hd}.1).

The second assertion (\ref{3.phenomena.hd}.2)
directly follows from Theorem \ref{codim.3.slc.GL.thm}.
 \qed
\end{say}

\section{Stability of exceptional divisors}
\label{exdivdef.sec}

We consider part 3 of  Theorem \ref{3.phenomena.hd}. Let
 $g:X\to C$ be a flat morphism to a smooth curve.
 Let $0\in C$ be a closed  point 
such that $X_0$ is $\q$-factorial. Let
  $Z_0\subset X_0$ be a subscheme of codimension $\geq 3$ and
$f:Y\to X$ be a projective, birational morphism
such that  $f$ is an isomorphism over $X\setminus Z_0$
and $f_0:Y_0\to X_0$ is birational but not an isomorphism.

Let $H_0\subset Y_0$ be an ample divisor.
Since $X_0$ is $\q$-factorial,  $m\cdot f_0(H_0)$ is Cartier
for some $m>0$. Thus
$$
E_0:=f_0^*\bigl(m f_0(H_0)\bigr)-mH_0
$$
is a nonzero,  $f_0$-exceptional Cartier divisor.
We will show that this implies that $Y_t\to X_t$ is not an isomorphism,
contrary to our assumptions.

 More generally, let $Y_0$ be a complex analytic space and
$E_0\subset Y_0$ a proper, complex analytic subspace.
We would like to prove that, under certain conditions, every deformation
$\{Y_t:t\in \dd\}$ induces a corresponding deformation
$\{E_t\subset Y_t:t\in \dd\}$.

If $E_0$ is a Cartier divisor, then by deformation theory
(see, for instance, \cite[Sec.I.2]{rc-book} or \cite[Sec.6]{har-def}) the
 obstruction space is $H^1\bigl(E_0, \o_{Y_0}(E_0)|_{E_0}\bigr)$.
If $E_0$ is smooth and its normal bundle is negative,
then by Kodaira's vanishing theorem the obstruction group is zero,
hence every flat deformation of $Y_0$
 induces a  flat deformation of the pair
$(E_0\subset Y_0)$.

Here we address the more general case when there is a  projective morphism
$f_0:Y_0\to X'_0$ which contracts $E_0$ to a point.
(This always  holds if the normal bundle of $E_0$ is negative, 
at least for analytic or algebraic spaces, see \cite{artin}.)
We allow $E_0$ to be singular.
By \cite{MR0304703},  any flat deformation
$\{Y_t:t\in \dd\}$ induces a corresponding deformation
$\{f_t:Y_t\to X_t:t\in \dd\}$ (with the slight caveat that
$X_0$ need not be normal, but its normalization is $X'_0$).
We can state our result in a more general form as follows.

\begin{prop}\label{stab.exc.divs.cor}
Let $g:X\to \dd$ be a flat morphism of pure relative dimension $n$.
Let $f:Y\to X$ be a projective, bimeromorphic morphism
such that $f_0:Y_0\to X_0$ is also bimeromorphic.

Assume that there is a nonzero (but not necessarily effective)
$\q$-Cartier divisor $E_0\subset Y_0$ such that
$\dim f_0\bigl(\supp E_0\bigr)\leq n-3$.

Then, for every $|t|\ll 1$ there is a   nonzero exceptional divisor 
$E_t\subset \ex(f_t)$.
\end{prop}

Proof. By taking general hyperplane sections of $X$ we may assume that
$f_0\bigl(\supp E_0\bigr)$ is a point $p\in X_0$. 

We use the notation  of Theorem \ref{G-L.setup.say}.
 Lemma \ref{G-L.setup.nonvanish.lem} shows that the cup product pairing
$$
H^{2}\bigl( V_0, \partial  V_0, \q\bigr)\times
H^{2n-2}\bigl( V_0,  \q\bigr)\to 
H^{2n}\bigl( V_0, \partial  V_0, \q\bigr)\cong \q
$$
is nonzero. On the other hand, if $\dim \ex(f_t)\leq n-2$ for $t\neq 0$
then Lemma \ref{G-L.setup.vanish.lem2} applies and so
Lemma \ref{G-L.setup.vanish.lem} shows that the above
cup product pairing is zero, a contradiction.\qed

\begin{rem}\label{divcont.def.applic.rem}
 (1) Note first that we do not assert that $\{E_t:t\in \dd\}$ is a
flat family of divisors, nor do we claim that the $E_t$ are $\q$-Cartier.
Most likely both of these hold under some natural hypotheses.

(2) The dimension restriction $\dim f_0\bigl(\supp E_0\bigr)\leq n-3$
is indeed necessary. If $Y_0$ is a smooth surface and
$E_0\subset Y_0$ is a smooth rational curve then the analog of 
Proposition \ref{stab.exc.divs.cor} holds only if $E_0$ is a $(-1)$-curve.

(3) The existence of a $\q$-Cartier divisor $E_0\subset Y_0$
seems an unusual assumption, but it is necessary, as shown by
the following examples.

Let $W$ be any smooth projective variety of dimension $n$
and $L$ a very ample line bundle
on $W$. Let $Y$ be the total space of the rank $r\geq 2$ bundle
$L^{-1}+\cdots +L^{-1}$ with zero section $W\cong E\subset Y$. 
Let $f:Y\to X$ be the contraction of $W$ to a point;
that is, $X$ is the spectrum of the symmetric algebra
$ H^0\bigl(W, \sym\bigl(L+\cdots +L\bigr)\bigr)$.

For any general map $g:X\to \a^1$ the conclusion of Proposition 
\ref{stab.exc.divs.cor} fails since we have $Y\to X$.
The fiber over the origin is a  hypersurface $Y_0\subset Y$
that contains $E$ and the codimension of $E$ in $Y_0$ is $r-1$.

If $r>n$ then a general $Y_0$ is smooth but for these
$\dim E\leq \frac12 \dim Y_0$.

If $r\leq n$ then $Y_0$ is always singular.
The most interesting case is when $r=n$. Then, 
for general choices, the only singularities of $Y_0$
are  ordinary nodes along $E$. 
If $n=r=2$ then $E$ is a divisor in $Y_0$ but it is not $\q$-Cartier
at these nodes.

An interesting special case arises when $W=\p^n$, $L=\o_{\p^n}(1)$ and
$r=n+1$. Then $X_0$ has a terminal singularity and
$Y_0\to X_0$ is crepant.

(4) The conclusion of Proposition \ref{stab.exc.divs.cor} should hold if $f$
is only proper, but the current proof uses projectivity in an essential way.

(5)  An examination of the proof shows that Proposition 
\ref{stab.exc.divs.cor}
can be extended to higher codimension exceptional sets as follows.
In view of the examples in (3), the assumptions seem to be optimal.
\end{rem}

\begin{prop} \label{divcont.def.applic.rem.6} 
Let $g:X\to \dd$ be a flat morphism of pure relative dimension $n$.
Let $f:Y\to X$ be a projective, bimeromorphic morphism
such that $f_0:Y_0\to X_0$ is also bimeromorphic.

Assume that  $\ex(f_0)$ is mapped to a point, $d:=\dim \ex(f_0)>n/2$ and 
$\ex(f_0)$ supports an effective $d$-cycle
that is the Poincar\'e dual
 of a  cohomology class in $H^{2(n-d)}(Y_0, \partial Y_0,\q)$.
(The latter always holds if $Y_0$ is smooth.)

Then $\dim \ex(f_t)=d$ for every $t$. \qed
\end{prop}

\section{Another topological approach}
\label{top-approach.sec}

The aim of this section is to recall the usual topological approach to
Problem \ref{main.g-l.conj}, going back at least to \cite{MR0215323}.
 This  works
if $X\setminus X_0$ is smooth. Then we discuss a possible
modification of the method that could lead to a proof over $\c$.

\begin{say}[Set-up]\label{leff.in.intersection.hom.setup}
 Let $X\subset \c^N$ be an analytic space of pure dimension $n$ and
$X_0\subset X$ a Cartier divisor.
Let $p\in X_0$ be a point,
 $W\subset X$  the intersection of $X$ with a small closed ball
around $p$ 
and set $W_0:=W\cap X_0$. We assume that
\begin{enumerate}
\item the interior of $W$ is Stein,
\item $W\setminus \{p\}$ is  homeomorphic to $\partial W\times [0,1)$,
\item $W_0\setminus \{p\}$ is  
homeomorphic to $\partial W_0\times [0,1)$ and
\item $X\setminus X_0$ is smooth.
\end{enumerate}

\end{say}

\begin{prop} \label{leff.in.intersection.hom}
Notation and assumptions as above. Then the natural map
$$
H^i(W  \setminus \{p\}, \z)\to H^i(W_0  \setminus \{p\}, \z)
$$
is an isomorphism for $i\leq n-3$ and an injection for
$i=n-2$. 
\end{prop}

Proof. By assumption, 
$$
H^i(W  \setminus \{p\}, \z)\cong H^i(\partial W , \z)
\qtq{and}
H^i(W_0  \setminus \{p\}, \z)\cong H^i(\partial W_0 , \z).
$$
We have an exact sequence
$$
H^i(\partial W , \partial W_0 ,\z)\to 
H^i(\partial W , \z)\to 
H^i(\partial W_0 , \z)\to 
H^{i+1}(\partial W , \partial W_0 ,\z)
\eqno{(\ref{leff.in.intersection.hom}.1)}
$$
Since $W\setminus W_0$ is smooth,  Poincar\'e duality shows that
$$
H^i(\partial W , \partial W_0 ,\z)=
H_{2n-1-i} (\partial W \setminus \partial W_0 ,\z)\cong
H_{2n-1-i} ( W \setminus  W_0 ,\z).
\eqno{(\ref{leff.in.intersection.hom}.2)}
$$
By assumption the interior of $W$ is $n$-dimensional and  Stein, hence
 so is  the interior of $ W \setminus  W_0$.
Thus $ H_{2n-1-i} ( W \setminus  W_0 ,\z)=0$ for 
$2n-1-i\geq n+1$ by Theorem \ref{stein.top.main.thm}.

If $i\leq n-3$ then both groups at the end of (\ref{leff.in.intersection.hom}.1)
are zero, giving the isomorphism 
$H^i(W  \setminus \{p\}, \z)\cong H^i(W_0  \setminus \{p\}, \z)$.

If  $i= n-2$ then only the group on the left vanishes, thus we get
only an injection 
$H^{n-2}(W  \setminus \{p\}, \z)\into H^{n-2}(W_0  \setminus \{p\}, \z)$.
 \qed

\begin{cor} Notation and assumptions as above.
Assume in addition that $\depth_p\o_{X_0}\geq 2$ and $\dim X_0\geq 3$. Then
the natural restriction
$$
\pic^{an}(W  \setminus \{p\})\to \pic^{an}(W_0  \setminus \{p\})
$$
is an injection. 
\end{cor}

Proof. Consider the commutative diagram
$$
\begin{array}{ccc}
\pic^{an}(W  \setminus \{p\})& \into  & H^2(W  \setminus \{p\}, \z)\\
\downarrow && \downarrow\\
\pic^{an}(W_0  \setminus \{p\})& \to  & H^2(W_0  \setminus \{p\}, \z)
\end{array}
$$
The top horizontal arrow is injective by 
Lemma \ref{pic.H2.lem} and
the right hand vertical arrow is injective by
Proposition \ref{leff.in.intersection.hom}, hence the left 
 hand vertical arrow is also  injective.\qed

\medskip

The next lemma shows that, even in the classical setting,
that is when $Y$ is projective and $Y_0\subset Y$ is an ample divisor,
the restriction map
$H^2(Y, \q)\to H^2(Y_0, \q)$ is not  injective under some conditions.
 We then show in 
(\ref{Lef.fails.for.H2.exmp}) that such examples do exist.

\begin{lem} \label{Lef.fails.for.H2.lem}
 Let $X$ be a smooth projective variety of dimension $n$ and
$Z\subset X$ a smooth divisor.
Assume that $H^1(X,\q)=0$ and there is a morphism
$g:X\to Y$ that  contracts $Z$ to a point
$y\in Y$ and  is an isomorphism otherwise.
Let $Y_0\subset Y$ be a smooth divisor (not passing through $y$). Then
the kernel of the restriction map $H^2(Y, \q)\to H^2(Y_0, \q)$
 contains $H^1(Z, \q) $.
\end{lem}

Proof. The cohomology sequence of the pair  $(Y, y)$ shows that
$H^2(Y, \q)\cong H^2(Y,y, \q)$ and $H^2(Y,y, \q)\cong H^2(X,Z, \q)$.
The latter in turn sits in an exact sequence
$$
H^1(X,\q)\to H^1(Z, \q)\to H^2(X,Z, \q)\to H^2(X, \q).
$$
We assumed that $H^1(X,\q)=0$ hence there is an injection
$$
 H^1(Z, \q)\into \ker\bigl[H^2(X,Z, \q)\to H^2(X, \q)\bigr].
$$
Since $H^2(Y, \q)\to H^2(Y_0, \q)$ factors as 
$$
H^2(Y, \q)\cong H^2(X,Z, \q)\to H^2(X, \q)\to H^2(Y_0, \q)
$$
we obtain an injection
$$
 H^1(Z, \q)\into \ker\bigl[H^2(Y, \q)\to H^2(Y_0, \q)\bigr].\qed
$$

We thus need to find examples as above where $H^1(Z, \q)\neq 0$.
In the next examples, $Z$ is an Abelian variety.

\begin{exmp}\label{Lef.fails.for.H2.exmp}
 Let $\pi_0: S\to \p^1$ be a simply connected  elliptic surface.
For instance, we can take $S$ to be the blow-up of $\p^2$
at the 9 base points of a cubic pencil.

By composing $\pi$ with suitable automorphisms of $\p^1$
we get $n$ simply connected  elliptic surfaces
$\pi_i:S_i\to \p^1$ such that for every point $p\in \p^1$
at most one of the $\pi_i$ has a singular fiber over $p$.
Thus the fiber product
$$
X_1:=S_1\times_{\p^1}S_2\times_{\p^1}\cdots \times_{\p^1}S_n
$$
is a smooth variety of dimension $n+1$. General fibers of the projection
$\pi_1:X_1\to \p^1$  are Abelian varieties of dimension $n$ and
$X_1$ is simply connected.

Fix  Abelian fibers $A_1, A_2\subset X_1$. 
Let $H_1\subset X_1$ be a very ample
divisor such that $A_1\cap H_1$ is smooth. 
Let $X:=B_{A_1\cap H_1}X_1\to X_1$ denote the blow-up of $A_1\cap H_1$.
Let $H\subset X$ denote the birational transform of $H_1$
and $A\subset X$  the birational transform of $A_1$.

For $m\gg 1$, the linear system  $|H+mA_2|$ is base point free.
This gives a morphism
$g:X\to Y$. Note that $g$ contracts $A$ to a point
$y\in Y$ and  $g:X\setminus A\to Y\setminus \{y\}$ is an isomorphism.
\end{exmp}


\def\cprime{$'$} \def\cprime{$'$} \def\cprime{$'$} \def\cprime{$'$}
  \def\cprime{$'$} \def\cprime{$'$} \def\dbar{\leavevmode\hbox to
  0pt{\hskip.2ex \accent"16\hss}d} \def\cprime{$'$} \def\cprime{$'$}
  \def\polhk#1{\setbox0=\hbox{#1}{\ooalign{\hidewidth
  \lower1.5ex\hbox{`}\hidewidth\crcr\unhbox0}}} \def\cprime{$'$}
  \def\cprime{$'$} \def\cprime{$'$} \def\cprime{$'$}
  \def\polhk#1{\setbox0=\hbox{#1}{\ooalign{\hidewidth
  \lower1.5ex\hbox{`}\hidewidth\crcr\unhbox0}}} \def\cdprime{$''$}
  \def\cprime{$'$} \def\cprime{$'$} \def\cprime{$'$} \def\cprime{$'$}
\providecommand{\bysame}{\leavevmode\hbox to3em{\hrulefill}\thinspace}
\providecommand{\MR}{\relax\ifhmode\unskip\space\fi MR }
\providecommand{\MRhref}[2]{%
  \href{http://www.ams.org/mathscinet-getitem?mr=#1}{#2}
}
\providecommand{\href}[2]{#2}

\bigskip

\noindent Princeton University, Princeton NJ 08544-1000

{\begin{verbatim}kollar@math.princeton.edu\end{verbatim}}

\end{document}